\documentclass{svjour2}                    
\smartqed
\usepackage[OT1]{fontenc}
\usepackage[colorlinks,citecolor=blue,urlcolor=blue]{hyperref}
\usepackage{latexsym}
\usepackage{graphicx}
\usepackage{float}
\usepackage{wrapfig}
\usepackage{amsmath}
\usepackage{amsfonts}
\usepackage{amssymb}

\newcommand{\detg}{\textrm{sdet}}
\newcommand{\tr}{\,\textrm{tr}\,}
\newcommand{\trg}{\,\textrm{str}\,}
\newcommand{\Det}{\textrm{det}}
\newcommand{\im}{\textrm{Im}}
\newcommand{\re}{\textrm{Re}}
\newcommand{\diag}{\textrm{diag}}

\newcommand{\eins}{\leavevmode\hbox{\small1\kern-3.8pt\normalsize1}}
\begin{document}

\title{Supersymmetry approach to Wishart correlation matrices: Exact results}

\titlerunning{Supersymmetry approach to Wishart correlation matrices}

\author{Christian Recher \and Mario Kieburg \and Thomas Guhr \and Martin R. Zirnbauer}

\authorrunning{C. Recher et al.} 

\institute{C. Recher \at
              Instituto de Ciencia de Materiales de Madrid  \\
              CSIC, Cantoblanco Madrid\\
              28049 Madrid, Spain\\
              \email{christian.recher@uni-due.de}
           \and
           M. Kieburg \at
              Department of Physics and Astronomy\\
              State University of New York \\
              Stony Brook, NY 11794-3800, USA\\
              \email{mario.kieburg@stonybrook.edu}
              \and
           T. Guhr \at
              Fakult\"at f\"ur Physik,
              Universit\"at Duisburg-Essen \\
              Lotharstra\ss e 1\\
              47048 Duisburg, Germany\\
              \email{thomas.guhr@uni-due.de}
              \and
           M. R. Zirnbauer \at
              Institut f\"ur Theoretische Physik,
              Universit\"at zu K\"oln \\
              Z\"ulpicher Stra\ss e 77\\
              50937 K\"oln, Germany\\
              \email{zirnbauer@uni-koeln.de}
}

\date{Received: date / Accepted: date}

\maketitle

\begin{abstract}
We calculate the `one-point function', meaning the marginal probability density function for any single eigenvalue, of real and complex Wishart correlation matrices. No explicit expression had been obtained for the real case so far. We succeed in doing so by using supersymmetry techniques to express the one-point function of real Wishart correlation matrices as a twofold integral. The result can be viewed as a resummation of a series of Jack polynomials in a non-trivial case. We illustrate our formula by numerical simulations. We also rederive a known expression for the one-point function of complex Wishart correlation matrices.
\end{abstract}

\section{Introduction}

Complex systems of many different kinds are in the focus of modern research \cite{GMW98,Ver04,Haa01,Voi05}. Correlation matrices obtained from data sampling are a key tool to study such systems \cite{lal99,Spr08}. A standard approach in multivariate statistics is to use random matrix theory (RMT) to model the correlation matrices \cite{Mui82}. Adopting the framework of RMT, we will calculate the marginal probability density function for any single eigenvalue to take a given value $x$ (referred to as the one-point function $S_\beta(x)$ for short) for the case of real ($\beta = 1$) and complex ($\beta = 2$) correlation matrices. In the complex case, an explicit expression is known \cite{Alf04}. Real correlation matrices are encountered more frequently, but unfortunately, closed expressions for their one-point functions and related quantities have not been obtained so far. This is because a certain integral over the orthogonal group is not available in explicit form. Sophisticated power series techniques have been developed in order to tackle this problem \cite{Mui82}. However, the resulting expressions suffer from the drawback that a resummation of the infinite series has not been possible so far. For correlation matrices of large dimension, asymptotic results were derived in \cite{Sil95}. Recently some new results for the one-point function and the two-point correlation function in the asymptotic regime have been found \cite{VP2010}.

Here we provide exact results for the one-point functions of real and complex correlation matrices. We use an alternative approach to  circumvent the problems mentioned above. Our approach relies on the supersymmetry method \cite{Efe97} -- nowadays a standard tool for RMT applications in physics \cite{Verb04}. We derive an exact expression for the one-point function of real correlation matrices as a twofold integral. We also rederive the known result for the one-point function of complex correlation matrices. In Ref.\ \cite{RKG10/1}, we presented our main results to make them available for applications. Here, we give the full derivation of our results in a form that addresses not only physicists and practitioners, but also the mathematics and statistics community.

The article is organized as follows. We formulate the problem and
introduce our notation in Sec.\ \ref{sec.2}.  In Sec.\ \ref{sec.3} we
apply the supersymmetry technique, pursuing in parallel two different
approaches put forward in the literature, namely the generalized Hubbard-Stratonovich transformation \cite{Guh06,KGG08} and the superbosonization formula \cite{Som07,LSZ07}. These two approaches are equivalent \cite{KSG08}, but they yield different forms of the final expressions, each having their own advantages and disadvantages. We explicitly calculate the one-point function in Sec.\ \ref{sec.4}. In Sec.\ \ref{sec.5} we numerically integrate our formula to compare with Monte-Carlo simulations. We summarize and conclude in Sec.\ \ref{sec.6}.

\section{Formulating the Problem}\label{sec.2}

In Sec.\ \ref{sec2.1} we define the ensemble of random matrices to be considered in this paper. In Sec.\ \ref{sec2.2} we introduce a generating function for the one-point function. This generating function will serve as the starting point for the supersymmetry method.

\subsection{Ensemble of Wishart correlation matrices and one-point function}\label{sec2.1}

We briefly sketch the RMT approach to correlation matrices as set up in Ref.\ \cite{Mui82}. We consider real and complex Wishart correlation matrices. The building block for these are rectangular $p \times n$ matrices which we denote by $W= [W_{jk}], $ with $j=1, \ldots ,p$ and $k=1, \ldots ,n$. The $p$ rows can be viewed as the model time series of length $n$. The $p\times p$ matrix $WW^{\dagger}$ is the model correlation matrix, also referred to as the Wishart correlation matrix. We always assume $p \leq n $. The entries of $W$ are either real or complex random variables. These two cases are labeled by the Dyson index $\beta$ taking the value $\beta =1 $ for real entries ($W_{jk} \in \mathbb{R}$) and $\beta =2 $ for complex entries ($W_{jk} \in \mathbb{C}$). For the joint probability distribution of the entries of $W$ one chooses the multivariate Gaussian distribution
\begin{equation}
 P_{\beta}(W, C)= D_{\beta} \exp \left(-\frac{\beta}{2}
 \tr W^{\dagger}C^{-1} W \right),
 \label{2.1}
\end{equation}
where $C$ is the empirical, \textit{i.e.} given correlation matrix. By construction, the ensemble-averaged Wishart correlation matrix $WW^{\dagger}/n$ equals $C$, the empirical one. The Gaussian assumption \eqref{2.1} is justified in most if not all situations of interest. The full measure is $P_{\beta}(W, C)d[W]$ where
\begin{equation}
 d[W] = \begin{cases}
 \displaystyle\prod\limits_{j=1}^{p}\prod\limits_{k=1}^{n}dW_{jk} &\textrm{for}\ \beta =1,\\
 \displaystyle\prod\limits_{j=1}^{p}\prod\limits_{k=1}^{n}d\re W_{jk} d \im W_{jk} &\textrm{for}\ \beta =2 , \end{cases}\label{2.31}
\end{equation}
is the corresponding volume element. This measure fulfills the invariance condition
\begin{equation}
 P_{\beta}(W, C)d[W]= P_{\beta}(UW, UCU^\dagger)d[UW] \label{2.1a}
\end{equation}
for an arbitrary orthogonal ($\beta=1$) or unitary ($\beta=2$) $p \times p$ matrix $U$. Since the domain of $W$ ($\mathbb{R}^{p\times n}$ for $\beta=1$, $\mathbb{C}^{p\times n}$ for $\beta=2$) is invariant under the transformation $W \mapsto UW$, we may replace $C$ by the diagonal matrix of its eigenvalues as long as invariant quantities such as the one-point function (see below) are studied. Thus, we set $C \equiv \Lambda$ where $\Lambda = \diag (\Lambda_1, \ldots , \Lambda_p)$ is the diagonal matrix containing the eigenvalues. By the definition of $C$ as a correlation matrix we have $\Lambda_j >0$.  The constant $D_{\beta}$ in Eq.\ \eqref{2.1} ensures the normalization of $P_{\beta}(W, C)d[W]$ to unity and is given by
\begin{eqnarray}
 D_{\beta} = \big( (2 \pi / \beta)^p \, \Det \Lambda \big) ^{- n \beta/2}\label{2.2}.
\end{eqnarray}
The set of random matrices $WW^{\dagger}$ with the entries of $W$ distributed according to Eq.\ \eqref{2.1} is referred to as the ensemble of Wishart correlation matrices (sometimes also as the correlated Wishart ensemble). We mention in passing that for the choice $\Lambda=\eins_p$, where $\eins_p$ denotes the $p \times p$ unit matrix, the ensemble defined by Eq.\ \eqref{2.1} is equivalent to the so-called Gaussian chiral random matrix ensemble, which serves as a model for the universal eigenvalue statistics of the Dirac operator in  Quantum Chromodynamics \cite{Ver94/2}.

The one-point function for the eigenvalues $\lambda_j$ of $WW^{\dagger}$ is defined by
\begin{align}\label{2.3}
 S_{\beta}(x) &:= \int d[W] P_{\beta}(W,\Lambda) \frac{1}{p} \sum_{j=1}^p \delta(x-\lambda_j) \cr &= \frac{1}{\pi p} \lim_{\varepsilon \rightarrow 0+} \im \int d[W] P_{\beta}(W,\Lambda) \tr\frac{\eins_{p}}{(x-\mathrm{i} \varepsilon)\eins_{p} - WW^{\dagger}},
\end{align}
where in the second line we have passed to an expression involving the resolvent. By the definition \eqref{2.3}, the one-point function is a function of $x$ which parametrically depends on the empirical eigenvalues $\Lambda_1, \ldots ,\Lambda_p$. We drop the dependence on $\Lambda$ in writing $S_{\beta}(x)$.

Having defined the object of interest, we briefly comment on why it is difficult to handle the real case $\beta=1$ by the traditional techniques of multivariate analysis. In the standard approach to calculating the one-point function \eqref{2.3} one makes a singular-value decomposition
\begin{eqnarray}
 W=UwV,
\end{eqnarray}
where $U\in{\rm O}(p) , \ V \in {\rm O}(n)$  for $\beta=1$ and $U\in {\rm U}(p) ,\ V \in {\rm U}(n)$ for $\beta=2$. The $p \times n$ matrix $w$ contains the singular values $w_j \in \mathbb{R}$ ($j = 1, \ldots,p$) of $W$. For the matrix $WW^{\dagger}$ this decomposition yields
\begin{eqnarray}
 WW^{\dagger} = U w^2 U^\dagger \quad\textrm{with}
 \quad w^2 = \diag (w_1^2, \ldots ,w_p^2) .
\end{eqnarray}
While the substitution $W \mapsto UW$ leaves the resolvent in Eq.\ \eqref{2.3} invariant, the diagonalizing matrix $U$ does not drop out of the probability distribution function $P_{\beta}(W,\Lambda)$ containing the matrix $\Lambda \not= \eins_p$. The decomposition thus leads to the group integral
\begin{equation}
 \Phi_{\beta}(\Lambda,w^2) = \int \exp \Big(-\frac{\beta}{2} \tr
 U^{\dagger} \Lambda^{-1} U w^2 \Big) d\mu(U)\label{eq4}.
\end{equation}
For $\beta=2$ this is the celebrated Harish--Chandra-Itzykson-Zuber
integral \cite{Har58,ItzZub80} for which there exists an explicit expression. On the other hand, for $\beta=1$ no simple expression is known. The only expression available \cite{Mui82} is an infinite series in terms of zonal polynomials or equivalently, Jack polynomials, which in turn are only known recursively. A resummation of the infinite series has not been possible so far. In the present work we circumvent this problem by using a supersymmetry approach. This allows us to derive an expression for the one-point function as a twofold integral.

\subsection{Generating function}\label{sec2.2}

The starting point for the supersymmetry approach is the generating function
\begin{equation}\label{2.4}
 Z_{\beta}(x_0,x_1) = \int d[W] P_{\beta}(W,\Lambda) \frac{\det(x_1 \eins_p - WW^{\dagger} )}{\det( x_0 \eins_p - WW^{\dagger} )},
\end{equation}
where $x_0 , x_1$ are complex variables, $x_0 \notin \mathbb{R}_+$. The one-point function can be computed from it by taking a derivative:
\begin{equation}\label{2.5}
 S_{\beta} (x) = (2\pi\mathrm{i} p)^{-1} \frac{\partial}{\partial x_1} \bigg|_{x_1 = x} \lim_{\varepsilon \rightarrow 0+} \big( Z_{\beta}(x- \mathrm{i} \varepsilon,x_1) - Z_{\beta}(x + \mathrm{i} \varepsilon,x_1) \big) .
\end{equation}
Note that the generating function $Z_\beta(x_0,x_1)$ equals unity at $x_0 = x_1$. In the following we derive simple and computationally useful expressions for it by applying the supersymmetry technique.

\section{Passing to superspace}\label{sec.3}

There are several ways to express the generating function as an integral over a suitable superspace. Of particular prominence are the generalized Hubbard-Stratono\-vich transformation put forward in Ref.\ \cite{Guh06,KGG08} and the superbosonization formula derived in Ref.\ \cite{LSZ07}. Superbosonization \cite{EST04} was first proposed in a field theoretical context. It was then explored how the supersymmetry method can be extended to arbitrary invariant random matrix ensembles. In the
unitary case, this problem was solved \cite{Guh06} by introducing the
generalized Hubbard-Stratono\-vich transformation. In Ref.\ \cite{LSZ07}, rigorous superbosonization was developed for all classical Lie symmetries (unitary, orthogonal, symplectic). The approach of Ref.\ \cite{Guh06} was then completed by transcribing it to the orthogonal and symplectic cases \cite{KGG08}. The equivalence of the superbosonization of Ref.\ \cite{LSZ07} and the generalized Hub\-bard-Stratonovich transformation of Refs.\ \cite{Guh06,KGG08} was demonstrated in \cite{KSG08}.

In Sec.\ \ref{sec3.1}, we write the ratio of determinants in the generating function as a Gaussian integral over a rectangular supermatrix. Then we carry out the ensemble average. To the reader not experienced with anticommuting variables, we recommend the introductory parts of Refs.\ \cite{Efe97,VWZ85,Guh10} and the book by Berezin \cite{Ber87}. In Sec.\ \ref{sec3.2}, we use a duality between dyadic ordinary matrices and dyadic supermatrices to express the result of the ensemble average as a supermatrix integral. After analyzing certain symmetries of this dyadic supermatrix we replace it by a supermatrix of the same symmetries but independent matrix elements by means of the generalized Hubbard-Stratonovich transformation (Sec.\ \ref{sec3.3}), and alternatively with the help of the superbosonization formula (Sec.\ \ref{sec3.4}).

\subsection{Ensemble average}\label{sec3.1}

The determinant in the denominator of Eq.\ \eqref{2.4} can be expressed
as a Gaussian integral over a vector comprising ordinary commuting
variables. The determinant in the numerator can be expressed as a
Gaussian integral over a vector with anticommuting entries
\cite{VWZ85,Efe97}. By combining both expressions we obtain a
representation for the ratio of determinants in Eq.\ \eqref{2.4} as a Gaussian integral over a rectangular supermatrix $A$ comprising both vectors:
\begin{equation}\label{3.1}
 \frac{\det\left(x_1 \eins_p -WW^{\dagger}\right)}{\det\left(x_0 \eins_p -WW^{\dagger}\right)} = \int d[A] \exp \left( \frac{\mathrm{i}\beta}{2} \trg (X A^\dagger A - A^\dagger WW^{\dagger} A) \right),
\end{equation}
where the diagonal matrix $X$ is given by $X = \mathrm{diag}(x_0,x_0, x_1,x_1)$ for $\beta = 1$ and $X = \mathrm{diag}(x_0,x_1)$ for $\beta = 2$. We take $\im\, x_0 > 0$ in order for the Gaussian integral to converge. The rectangular supermatrix
\begin{eqnarray}\label{3.2.a}
\begin{array}{rclrcll}
 A & = & \begin{bmatrix} u_a, & v_a, & \zeta_a^*, & \zeta_a \end{bmatrix}_{1\leq a\leq p}, &\quad A^\dagger & = & \begin{bmatrix} u_b \\ v_b \\ \zeta_b \\ -\zeta_b^* \end{bmatrix}_{1\leq b\leq p} &\textrm{for}\  \beta=1,\\ A & = & \begin{bmatrix} z_a^*  &\zeta_a^* \end{bmatrix}_{1\leq a\leq p}, & A^\dagger & = & \begin{bmatrix} z_b \\ \zeta_b \end{bmatrix}_{1\leq b\leq p} &\textrm{for}\ \beta  =2,
\end{array}
\end{eqnarray}
is $p\times(2/\beta|2/\beta)$ dimensional. Here $u_j, v_j \in \mathbb{R}$ and $z_j \in \mathbb{C}$ are ordinary real or complex variables while $\zeta_j, \zeta_j^*$ are anticommuting variables, also referred to as Grassmann variables. We denote by $z_j^*$ the complex conjugate of $z_j$. In Eq.~\eqref{3.1}, $d[A]$ denotes the following product:
\begin{equation}\label{3.2.d}
 d[A]=\begin{cases} (2\pi)^{-p} \prod\limits_{j=1}^{p} du_j \, dv_j\, \partial_{\zeta_j^*} \partial_{\zeta_j} &{\rm for}\ \beta=1\,,\\
 \pi^{-p} \prod\limits_{j=1}^{p} d\re\, z_j\, d\im\, z_j\, \partial_{\zeta_j^*} \partial_{\zeta_j} &{\rm for}\ \beta=2\,. \end{cases}
\end{equation}
It should be mentioned that, in this context, one often writes $d\zeta \equiv \partial_\zeta$ for a Grassmann variable $\zeta$. However, the transformation law $d (t \zeta) = t^{-1} d\zeta$ for $t \in \mathbb{C}$ shows that $d\zeta = \partial / \partial\zeta$ really is a partial derivative, \emph{not} a differential!

By inserting the representation \eqref{3.1} into the generating function \eqref{2.4} and changing the order of doing the integrals we find
\begin{eqnarray}
 Z_{\beta}(x_0,x_1) &=& \int d[A] \exp \left( \frac{\mathrm{i}\beta}{2} \trg X A^\dagger A \right)\nonumber \\ &&\times D_{\beta} \int d[W] \exp\left( -\frac{\beta}{2} \tr W^{\dagger} (\Lambda^{-1}+\mathrm{i} AA^\dagger )W\right)\nonumber\\ &=& \int d[A] \exp \left( \frac{ \mathrm{i}\beta}{2} \trg X A^\dagger A \right) {\det}^{-n\beta/2}\left( \eins_p + \mathrm{i} AA^\dagger \Lambda \right).\label{3.2}
\end{eqnarray}
In the last step, we performed the Gaussian integral over $W$.

\subsection{Duality between ordinary and superspace}\label{sec3.2}

We now rewrite the determinant in Eq.\ \eqref{3.2} as a superdeterminant. This is possible due to duality relations between ordinary spaces and superspaces, see Refs.\ \cite{Zir96,Guh06,KGG08}. In the present context the duality amounts to the relation \cite{Guh06,KGG08}
\begin{eqnarray}\label{3.6}
 \det\left( \eins_p + \mathrm{i} AA^{\dagger} \Lambda \right)= \detg \left( \eins_{4/\beta} + \mathrm{i} A^{\dagger}\Lambda A\right).
\end{eqnarray}
We notice that the determinant is a polynomial while the superdeterminant is in principle a rational function. The relation originates from the identity $\tr (A A^\dagger \Lambda) = \trg (A^\dagger \Lambda A)$ and a Taylor expansion in the Grassmann variables, which is always a finite sum. The supermatrix $A^\dagger \Lambda A$ has dimension $4\times 4$ and $2\times 2$ for $\beta=1$ and $\beta=2$, respectively. On the other hand, the original matrix $AA^\dagger$ is $p\times p$ dimensional. This dimensional reduction is the crucial advantage of the supersymmetry method.

For present use, we take a look at the symmetry properties of the supermatrix $A^\dagger \Lambda A$. We see that in both cases ($\beta = 1, 2$) the left upper block (a.k.a.\ the boson-boson block) of $ A^{\dagger} \Lambda A$ is a Hermitian matrix. This observation will constrain some of the matrix blocks appearing below. What about the complex linear symmetries (\textit{i.e.} those not involving complex conjugation)? For $\beta = 2$ there are no such symmetries, but for $\beta = 1$ we have
\begin{eqnarray}\label{3.9}
 (A^{\dagger}\Lambda A)^T = S^T A^{\dagger}\Lambda A S \quad \textrm{with}\quad S = 
\begin{bmatrix} 1 & 0 & 0 & 0 \\ 0 & 1 & 0 & 0 \\ 0 & 0 & 0 & -1 \\ 0 & 0 & 1 & 0 \\ \end{bmatrix} \qquad ,
\end{eqnarray}
reflecting the fact that the related $p\times p$ matrix $\Lambda^{1/2} AA^{\dagger}\Lambda^{1/2}$ is symmetric.

Our aim now is to replace the supermatrix $ A^{\dagger}\Lambda A$
by a supermatrix $\sigma$ with independent matrix elements. We have
two approaches at our disposal: the generalized
Hubbard-Stratonovich transformation \cite{Guh06,KGG08}, and the superbosonization formula derived in Refs.\ \cite{Som07,LSZ07}.

\subsection{Generalized Hubbard-Stratonovich transformation} \label{sec3.3}

We proceed by introducing a (super-)Fourier representation of the required power of the superdeterminant function on the right-hand side of Eq.\ \eqref{3.6}: \begin{equation}\label{3.10}
 \detg^{-n\beta/2}(\eins_{4/\beta} +\mathrm{i} A^{\dagger}\Lambda A) = \int d[\varrho] I_\beta(\varrho) \exp\left(- \frac{\mathrm{i}\beta}{2} \trg A^{\dagger}\Lambda A\varrho \right) .
\end{equation}
The Fourier transform $I_\beta(\varrho)$ is
\begin{equation}\label{3.10-mrz}
 I_\beta(\varrho) = \int d[\sigma] \,\detg^{-n\beta/2}(\eins_{4/\beta} + \mathrm{i}\sigma)\exp\left(\frac{\mathrm{i}\beta}{2} \trg \sigma\varrho\right).
\end{equation}
In order for this integral representation to be formally consistent, the two supermatrices $\sigma$ and $\varrho$ have to share the complex linear symmetries of $A^{\dagger}\Lambda A$. Hence $\sigma$ and $\varrho$ are supermatrices of dimension 4$\times$4 for $\beta=1$ and 2$\times$2 for $\beta=2$, and for $\beta=1$ the complex linear constraint \eqref{3.9} is imposed.

We write the supermatrix $\sigma$ as
\begin{eqnarray}
 \sigma= \begin{bmatrix} \sigma_0 &\chi \\ \tilde\chi &~\mathrm{i} \sigma_1 \eins_{2/\beta} \end{bmatrix}, \label{3.11}
\end{eqnarray}
where the entries are $2/\beta$ dimensional square matrices. The matrix $\sigma_0$ (akin to the boson-boson block of $A^\dagger \Lambda A$) is Hermitian, and $\sigma_1$ is a real scalar. (The reason for putting the imaginary unit in front of $\sigma_1$ will become clear presently.) The off-diagonal blocks $\chi$ and $\tilde\chi$ contain all anticommuting variables of $\sigma$. For $\beta=2$ the diagonal block $\sigma_0$ is simply a real number whereas $\chi$ and $\tilde\chi \equiv \chi^*$ are two Grassmann variables. For $\beta=1$ the diagonal blocks $\sigma_0$ and $\sigma_1 \eins_2$ are real symmetric $2\times 2$ matrices. The off-diagonal blocks for $\beta =1$ have the following structure:
\begin{eqnarray} \label{3.11a}
 \chi =\begin{bmatrix}
       \eta & \eta^*\\ \xi & \xi^* \end{bmatrix},\qquad
 \tilde\chi = \begin{bmatrix}
       \eta^* & \xi^*\\ -\eta & -\xi \end{bmatrix},
\end{eqnarray}
where $\eta$ and $\xi$ denote Grassmann variables. With this choice $\sigma$ satisfies the constraint \eqref{3.9}. The supermatrix $\varrho$, similar to $\sigma$, is chosen as
\begin{eqnarray} \label{3.11b}
 \varrho = \begin{bmatrix} \varrho_0 &\omega \\
 \tilde\omega & ~\mathrm{i} \varrho_1 \eins_{2/\beta} \end{bmatrix}.
\end{eqnarray}
In a self-evident way, $\varrho$ is divided into blocks having the same symmetries as the blocks of $\sigma$ in Eq.\ \eqref{3.11}.

The super-integration measure for $\sigma$ in Eq.\ \eqref{3.10} is flat and reads
\begin{eqnarray}\label{eq3.24}
 d[\sigma]=\begin{cases} (2\pi)^{-2} d\sigma_{0aa}\, d\sigma_{0ab}\, d\sigma_{0bb}\, d\sigma_{1}\, \partial_{\eta} \partial_{\eta^*} \partial_{\xi} \partial_{\xi^*} & {\rm for} \quad \beta=1,\\ (2\pi)^{-1} d\sigma_0\, d\sigma_1\, \partial_\chi \partial_{\chi^*} &{\rm for }\quad\beta=2, \end{cases}
\end{eqnarray}
where $\sigma_{0aa}$ and $\sigma_{0bb}$ are the diagonal elements and $\sigma_{0ab}$ is the off-diagonal element of the real symmetric matrix $\sigma_0$ for $\beta = 1$. The measure $d[\varrho]$ is defined by an identical expression.

We now insert the representation \eqref{3.10} into Eq. \eqref{3.2}. Our ensemble-averaged generating function then becomes a supermatrix integral:
\begin{eqnarray}
 Z_{\beta}(x_0,x_1) &=& \int d[\varrho] I_\beta(\rho) \int d[A] \exp \left( \frac{\mathrm{i}\beta}{2} \trg ( X A^\dagger A - A^\dagger \Lambda A\varrho ) \right) \nonumber \\ &=& \int d[\varrho] I_\beta(\varrho) \prod_{j=1}^p \detg^{-\beta/2}\left( X -\varrho \Lambda_j \right),\label{3.13}
\end{eqnarray}
where in the last step we performed the integrals over $A$. Eq.\ \eqref{3.13} is the desired superspace representation of the generating function. Originally, the generating function was an integral over ordinary $p\times n$ matrices $W$. The representation \eqref{3.13} is an integral over supermatrices $\varrho$ of dimension $4 \times $4 for $\beta=1$ and $ 2\times 2$ for $\beta=2$. This drastically reduces the number of integrals to be calculated.

To complete the description of our result \eqref{3.13}, we must discuss the function $I_\beta (\varrho)$. This is a super-version of what has come to be called the Ingham-Siegel integral \cite{Guh06,KGG08}. (The name is due to Fyodorov \cite{Fyo02} who introduced such an integral in a related, non-super context.) An important point to appreciate is that the superdeterminant under the integral sign of \eqref{3.10-mrz} depends \emph{polynomially} on the variable $\sigma_1$. Therefore the Fourier transform \eqref{3.10-mrz} does not converge in that variable and the supersymmetric Ingham-Siegel integral $I_\beta(\varrho)$ cannot exist as a regular function of $\varrho$. Nevertheless, we can make sense of $I_\beta(\varrho)$ as a distribution in $\varrho$, as follows. (In the next subsection, we will reproduce the same result \eqref{3.13} by employing convergent integrals only.)

The Ingham-Siegel integral is invariant under conjugation of $\sigma$ by elements of the supergroups ${\rm UOSp}(2|2)$ or ${\rm U}(1|1)$ for $\beta = 1, 2$, respectively. This is \emph{not} automatic, but does hold true once the integral is properly regularized by an \emph{invariant} cutoff function. Thus, choosing an invariant cutoff of Gaussian form, we define the supersymmetric Ingham-Siegel integral by
\begin{equation}\label{3.10a-mrz}
 I_\beta(\varrho) := \lim_{\varepsilon\to 0+} \int d[\sigma] \, \detg^{-n\beta/2}(\eins_{4/\beta} +\mathrm{i}\sigma)\exp \left( \frac{\mathrm{i}\beta}{2} \trg \sigma\varrho\right) \,\mathrm{e}^{- \varepsilon \trg \sigma^2} .
\end{equation}
Here we see the reason why the block $\sigma_1 \eins_2$ is multiplied by $\mathrm{i} = \sqrt{-1}$: this factor cancels the minus sign from the supertrace, thereby making the Gaussian cutoff function decrease with increasing real integration variable $\sigma_1\,$.

By construction, the invariantly regularized Fourier transform \eqref{3.10a-mrz} is invariant and thus the distribution $I_\beta( \varrho)$ depends only on the eigenvalues of $\varrho = UR\,U^{-1} = U\diag(R_0, R_1\eins_{2/ \beta}) U^{-1}$. We note that the diagonal matrix $R_0$ has the dimension $2/\beta$ while $R_1$ is a scalar. An explicit formula for the Ingham-Siegel distribution $I_\beta(\varrho)$ was derived in Refs.\ \cite{Guh06,KGG08}. The result is
\begin{align}\label{A.11}
 I_{\beta}(\varrho) &= K_{\beta}\,\Theta(R_0) \, \mathrm{det}^{(n+1) \beta/2 -1}(R_0)\cr &\times \exp\left(- \frac{\beta}{2} \trg R \right) \left( \mathrm{i} \frac{\partial}{\partial R_1} \right)^{n-2/\beta} \delta(R_1) ,
\end{align}
with the constants
\begin{eqnarray}\label{A.11a}
 K_1 = \frac{\pi}{(n-2)!}, \qquad K_2 = \frac{2\pi}{(n-1)!} ,
\end{eqnarray}
and the Heaviside distribution
\begin{eqnarray}
 \Theta(R_0) = \begin{cases} 1 &\textrm{if}\  R_0\ \textrm{is a positive definite matrix,}\\ 0 &\textrm{else} . \end{cases}
\end{eqnarray}
This completes our description of the result \eqref{3.13}.

\subsection{Approach using superbosonization}\label{sec3.4}

We now rederive the result \eqref{3.13} (or, more precisely, an equivalent formula) by a different approach, avoiding the use of super-distributions. This will be mathematically clean in every respect. The price to pay is that we rely on a ``black box'', namely the superbosonization formula proved in \cite{LSZ07}.

We start over from the very beginning, Eq.\ \eqref{3.1}. Motivated by the fact that our probability measure $P_\beta(W,\Lambda) d[W]$ for $W \in \mathbb{C}^{p\times n}$ (resp.\ $W \in \mathbb{R}^{p\times n}$) is right-invariant but not left-invariant, we first pass from determinants on the left space $\mathbb{C}^p$ ($\mathbb{R}^p$) to determinants on the right space $\mathbb{C}^n$ ($\mathbb{R}^n$), and only afterwards introduce the standard Gaussian integral representation:
\begin{align}
 &\left( \frac{x_0}{x_1} \right)^{p-n} \frac{\det\left( x_1 \eins_p - WW^{\dagger}\right)} {\det\left(x_0 \eins_p - WW^{\dagger} \right)} =
 \frac{\det\left( x_1 \eins_n - W^{\dagger} W\right)} {\det\left(x_0 \eins_n - W^{\dagger} W \right)} \cr &= \detg^{-\beta/2} \left(X \otimes \eins_n - \eins_{4/\beta} \otimes W^\dagger W \right) \cr &= \int d[A] \exp \left( \frac{\mathrm{i}\beta}{2}\trg (X A^\dagger A - A^\dagger W^{\dagger} W A) \right) \label{mrz-sb1}.
\end{align}
$X$ is the same diagonal matrix as before, and we still take $\im\, x_0 > 0$ in order for the Gaussian integral to converge. The rectangular supermatrix $A$ is the same as before except for the change in dimension $p \to n$.

Next we take the expectation with respect to the probability measure $P_\beta (W, \Lambda) d[W]$ of the $W$-dependent factor under the integral sign:
\begin{align}
 &\int d[W] P_\beta (W, \Lambda) \exp \left( -\frac{\mathrm{i}\beta}{2} \tr  W A A^\dagger W^\dagger \right) \cr &= {\det}^{-\beta/2}\left( \eins_n \otimes \eins_p + \mathrm{i} A A^\dagger \otimes \Lambda \right) \cr &= {\detg}^{- \beta/2}\left( \eins_{4/\beta} \otimes \eins_p + \mathrm{i} A^\dagger A \otimes \Lambda \right)
 \label{mrz-sb2}.
\end{align}
In the last step we invoked the duality of Sec.\ \ref{sec3.2}. For the generating function \eqref{2.4} we then get the formula
\begin{align}\label{mrz-sb3}
 Z_\beta(x_0,x_1) = (x_0/x_1)^{n-p} &\int d[A]\exp\left( \frac{ \mathrm{i} \beta}{2} \trg X A^\dagger A \right) \cr &\times \prod_{j=1}^p {\detg}^{ -\beta/2}\left(\eins_{4/\beta}+\mathrm{i} A^\dagger A \Lambda_j \right).
\end{align}

Now a beautiful feature of the integrand in \eqref{mrz-sb3} is that it depends only on the product $A^\dagger A$ which is invariant under left translations $A\mapsto UA$ by $U\in\mathrm{O}(n)$ (resp.\ $\mathrm{U}(n)$) for $\beta = 1$ ($\beta = 2$). This means that we are exactly in the situation where the superbosonization formula of \cite{Som07,LSZ07} applies.

In a nutshell, superbosonization lets us switch from our $A$-integral to the same integral over a supermatrix $Q$ replacing $A^\dagger A$. The result is
\begin{align}\label{mrz-sb4}
 &Z_\beta(x_0,x_1) = (x_0/x_1)^{n-p} \int DQ \, F_\beta(Q) , \\ &F_\beta(Q) = \detg^{n\beta/2}(Q) \exp\left( \frac{\mathrm{i}\beta} {2}\trg X Q\right) \prod_{j=1}^p {\detg}^{-\beta/2} \left(\eins_{4/\beta} + \mathrm{i} Q \Lambda_j \right). \nonumber
\end{align}
$Q$ is formally identical to the supermatrix $\varrho$ of Sec.\ \ref{sec3.3}. The only difference is in the domain of integration: the present $Q$-integral is over a Riemannian symmetric superspace (of Cartan type $A$) with invariant (or Berezin-Haar) measure $DQ$. This is to say that $Q$ runs over the positive matrices in the left upper (or boson-boson) block and the unitary matrices in the right lower (or fermion-fermion) block. For $\beta = 2$ the measure $DQ$ is the flat measure $d[Q] \equiv d[\sigma]$ of \eqref{eq3.24} (up to a normalization factor), for $\beta = 1$ it is the flat one times an extra factor of $\detg^{-1/2}(Q)$.

Our two expressions \eqref{3.13} and \eqref{mrz-sb4} become identical upon making the substitution $Q \to \mathrm{i} Q / X$. In fact, the Heaviside distribution in $R_0$ of \eqref{A.11} restricts the integration over all Hermitian matrices $\varrho_0$ to the positive ones, and the derivatives of the $\delta$-distribution in $R_1$ of \eqref{A.11} have exactly the same effect as doing the integral over the unitary variable in the fermion-fermion block of $Q$ by the residue theorem. Thus the $\varrho$-integral in \eqref{3.13} is effectively over the Riemannian symmetric superspace parameterized by $Q$; see \cite{KSG08} for more discussion of how the two approaches are related.

\section{Explicit expressions for the one-point function}\label{sec.4}

We consider in Sec.\ \ref{sec4.1} the complex case ($\beta=2$) and
rederive the result found in Ref.\ \cite{Alf04}. In Sec.\ \ref{sec4.2} we study the real case ($\beta = 1$) and derive an expression as a twofold integral. In both cases we have a choice between the methods of superbosonization and generalized Hubbard-Stratonovich transformation.

\subsection{Complex case}\label{sec4.1}

We show how to reproduce in a few more steps the known result \cite{Alf04} for the one-point function for $\beta = 2$. The supermatrix $Q$ in this case is of size $2 \times 2$. We begin by introducing its eigen-representation:
\begin{equation}\label{mrz-eigen}
 Q = U R\,U^{-1}, \quad U = \begin{bmatrix} \sqrt{1 - \alpha \alpha^\ast} &- \alpha\cr \alpha^\ast ~ &\sqrt{1 - \alpha^\ast \alpha}\end{bmatrix}, \quad R = \begin{bmatrix} r &0\cr 0 &s\end{bmatrix},
\end{equation}
with eigenvalues $r, s$ and two Grassmann variables $\alpha, \alpha^\ast$. This parametrization is singular (more precisely, degenerate in the Grassmann variables) at $r = s$. To suppress any effects of the singularity, we are going to utilize the scale invariance of $DQ$ to change the integration radius $|s| = 1$ to $|s| = q $ with $q \to \infty$. By the change-of-variables formula for superintegrals, our invariant integral \eqref{mrz-sb4} for the choice of parametrization \eqref{mrz-eigen} then takes the form
\begin{equation}\label{mrz-3.ii}
 \int DQ\, F_2(Q) = (2\pi\mathrm{i})^{-1} \lim_{q \to \infty} \int\limits_0^q \! dr \!\! \oint\limits_{|s|=q} \!\! ds \; (r-s)^{-2} \partial_\alpha \partial_{\alpha^\ast} F_2 \,,
\end{equation}
provided that $F_2(Q)$ vanishes on the locus of the coordinate singularity $r = s = q \to \infty$. (Otherwise, so-called Efetov-Wegner boundary terms appear.) To arrange for this vanishing property to hold, notice that for large values of $r = s = q$ the function $F_2(q \eins_2)$ behaves as $\mathrm{e}^{ \mathrm{i} q (x_0 - x_1)}$. If $\im (x_0 - x_1) > 0$ then this exponential factor makes $F_2(q \eins_2)$ vanish on the singular locus for $q \to +\infty$. We therefore assume the inequality $\im\,x_0 > \im\,x_1$.

Since the two superdeterminants in $F_2(Q)$ are functions of the eigenvalues $r$ and $s$ only, all of the dependence on $\alpha, \alpha^\ast$ resides in the factor $\mathrm{e}^{\mathrm{i} \trg X Q}$ and the process of integrating  (actually, differentiating) w.r.t.\ the Grassmann variables is simply done by
\begin{equation}\label{mrz-3.iii}
 \partial_\alpha \partial_{\alpha^\ast}\,\mathrm{e}^{\mathrm{i} \trg X Q} = \mathrm{e}^{\mathrm{i}(x_0 r - x_1 s)}\mathrm{i}(x_0 - x_1)(r-s).
\end{equation}
By using \eqref{mrz-3.ii} and \eqref{mrz-3.iii} in \eqref{mrz-sb4} we obtain
\begin{equation}\label{mrz-3.iv}
 Z_2(x_0,x_1) = (x_0 / x_1)^{n-p} \lim_{q\to\infty} \int\limits_0^q \! dr \!\! \oint\limits_{|s|=q}\!\!\! \frac{ds} {2\pi} \; \frac{r^n (x_0 - x_1) \, g_\Lambda(x_1;s)}{s^n (r-s) \, g_\Lambda(x_0;r)} ,
\end{equation}
where $g_\Lambda(x;s)$ is the function
\begin{equation}\label{def-gL}
 g_\Lambda(x;s) = \mathrm{e}^{-\mathrm{i}xs} \prod_{j=1}^p (1 + \mathrm{i}s \Lambda_j) .
\end{equation}

Our next step is to perform the inner integral over the circle variable $s$. This is done by invoking Cauchy's integral theorem to show that for any complex-analytic function $f(s)$ one has
\begin{equation}\label{eq:cauchy}
 \frac{1}{2\pi\mathrm{i}} \oint\limits_{|s| = q} \frac{f(s)\,ds}{s^n (r-s)} = \left\{ \begin{array}{ll} - r^{-n} f^{[\geq n]}(r), &\qquad q > r , \cr + r^{-n} f^{[< n]}(r), &\qquad q < r , \end{array} \right.
\end{equation}
where $f^{[\geq n]}(r)$ stands for the Taylor series of $f(r)$ in $r$ with all terms up to order $n-1$ deleted, and $f^{[<n]}(r) = f(r) - f^{[\geq n]}(r)$. Application of this formula to \eqref{mrz-3.iv} yields the result
\begin{equation}\label{mrz-3.vx}
 Z_2(x_0,x_1) = - \mathrm{i}\,(x_0 -x_1)(x_0 /x_1)^{n-p} \int \limits _0^\infty \! dr\,\frac{g_\Lambda^{[\geq n]}(x_1;r)} {g_\Lambda(x_0;r)} .
\end{equation}
Recall that $\im\, x_0 > \im\, x_1$ is required in order for this to hold. The correct expression for the opposite case of $\im\, x_0 < \im\, x_1$ is obtained by replacing in \eqref{mrz-3.vx} every occurrence of $\mathrm{i}$ by $- \mathrm{i}$ or, equivalently, by sending $r \to -r$.

We now turn to the calculation of the one-point function $S_2(x)$. We recall the expression \eqref{2.5} involving $Z_2 (x \pm \mathrm{i} \varepsilon, x_1)$ and apply our result \eqref{mrz-3.vx} to the case of $x_0 = x \pm \mathrm{i} \varepsilon$ with real parameters $x, x_1$ and $\varepsilon > 0$. Our integral representations of $Z_2(x \pm \mathrm{i} \varepsilon , x_1)$ then combine to a \emph{single} integral:
\begin{align}\label{eq:singleI}
 &\lim_{\varepsilon\to 0+} \big( Z_2(x - \mathrm{i} \varepsilon,x_1) - Z_2(x+\mathrm{i}\varepsilon,x_1) \big) \cr &= \mathrm{i}\,(x - x_1)(x / x_1)^{n-p} \lim_{\varepsilon\to 0+} \int \limits_{-\infty}^\infty \! dr\,\mathrm{e}^{-\varepsilon |r|}\, \frac{g_\Lambda^{[\geq n]}(x_1;r)} {g_\Lambda(x;r)} .
\end{align}
If we re-express the numerator of the integrand as $g_\Lambda^{[\geq n]} \equiv g_\Lambda - g_\Lambda^{[< n]}$ then the term $g_\Lambda$ contributes $(x-x_1) \delta(x-x_1) = 0$. Hence we replace this numerator by $- g_\Lambda^{[<n]} (x_1;r)$, which is a polynomial in $r$ of degree $n-1$.

Now for $x < 0$ we may close the integration contour around the lower half of the complex $r$-plane. From \eqref{def-gL} one sees that the integrand is holomorphic in $r$ for $\im\, r < 0$. We therefore get $S_2(x) = 0$ for $x < 0$. On the other hand, for $x > 0$ the contour has to be closed around the upper half-plane. Again, recall Eq.\ \eqref{2.5}. Because the $r$-integral is now manifestly finite for all values of $x_1$, we may safely take the $x_1$-derivative at $x_1 = x$ by simply removing the prefactor $x_1 - x$ and replacing $x_1$ by $x$ under the integral sign. (Please be warned that in \eqref{mrz-3.vx} the $r$-integral has to diverge at $x_1 = x_0$ to arrange for $Z_2(x_0,x_0) = 1$ in spite of the factor $x_0 - x_1$.) The result for $x > 0$ is
\begin{align}\label{mrz-40}
 S_2(x) = \frac{1}{2\pi p} \oint\! dr\, \frac{g_\Lambda^{[< n]}(x;r)} {g_\Lambda(x;r)} = - \frac{1}{p} \sum_{j=1}^p (\mathrm{i} r + 1/ \Lambda_j) \frac{g_\Lambda^{[\geq n]}(x;r)} {g_\Lambda(x;r)} \bigg\vert_{r \to \mathrm{i}/\Lambda_j} ,
\end{align}
where the integration contour in the first expression encloses the poles of $1 / g_\Lambda(x;r)$ at $r = \mathrm{i} / \Lambda_j$ $(j = 1, \ldots, p)$. To get the second expression, we switched from $g_\Lambda^{[< n]}$ back to $- g_\Lambda^{[\geq n]}$ and applied the residue theorem.

The residues are best computed by re-inserting the integral representation \eqref{eq:cauchy} for $g_\Lambda^{[\geq n]}$, which gives
\begin{equation}\label{mrz-3.vi}
 - g_\Lambda^{[\geq n]} (x;\mathrm{i}/\Lambda_j) = (\mathrm{i} / \Lambda_j)^{n-1} \oint \frac{ds}{2\pi\mathrm{i}} \, s^{-n} \mathrm{e}^{ -\mathrm{i}xs} \prod_{l (\not= j)} \left( 1 + \mathrm{i}s \Lambda_l \right) .
\end{equation}
We now use the expansion
\begin{equation}\label{mrz-3.v}
 \prod_{l (\not= j)} (1 + \mathrm{i} s \Lambda_l) = \sum_{k = 1}^p (\mathrm{i} s)^{k-1} E_{k-1}(\Lambda^{\widehat{j}})
\end{equation}
where $E_k$ are the elementary symmetric functions
\begin{equation}\label{eq:ESF}
 E_k(\Lambda) := \sum_{1 \leq j_1 < j_2 < \ldots < j_k \leq p} \Lambda_{j_1} \Lambda_{j_2} \cdots \Lambda_{j_k} \,.
\end{equation}
The notation $\Lambda^{\widehat{j}}$ means that the eigenvalue $\Lambda_j$ is to be dropped from the diagonal matrix $\Lambda$. Now, inserting the expansion \eqref{mrz-3.v} into \eqref{mrz-3.vi} we carry out the $s$-integral. Our final result for the one-point function then reads
\begin{equation}\label{mrz-final2}
 S_2(x) = \frac{\Theta(x)}{p} \sum_{j=1}^p \frac{\mathrm{e}^{-x/\Lambda_j} (1/\Lambda_j)^n} {\prod_{l (\not= j)} (1 - \Lambda_l / \Lambda_j)} \sum_{k=1}^p \frac{x^{n-k}}{(n-k)!}E_{k-1}(-\Lambda^{\widehat{j}}).
\end{equation}

To conclude this subsection, we present an alternative expression for the result \eqref{mrz-final2} as a ratio of determinants:
\begin{equation}\label{eq:mario}
 S_2(x) = \frac{\Theta(x)}{p} \det \begin{bmatrix} 0 &B_\Lambda(x) \cr C(x) & D_\Lambda \end{bmatrix} / \det D_\Lambda \;,
\end{equation}
where $D_\Lambda$ is the $p \times p$ matrix with matrix elements $(D_\Lambda)_{k,j} = \Lambda_j^{-k+1}$, and $B_\Lambda(x)$ and $C(x)$ are row and column vectors with entries
\begin{equation}
 B_\Lambda(x)_{j} = \mathrm{e}^{- x / \Lambda_j} (1/\Lambda_j)^n , \qquad C(x)_{k} = - \frac{x^{n-k}}{(n-k)!} .
\end{equation}
The denominator $\det D_\Lambda = \prod_{j > j^\prime} (\Lambda_j^{-1} - \Lambda_{j^\prime}^{-1})$ is essentially the Vandermonde determinant associated with the numbers $\Lambda_1 , \ldots, \Lambda_p$.

To verify the expression \eqref{eq:mario} one expands the $(p+1)\times (p+1)$ determinant with respect to the first row and first column and then uses a standard identity \cite{WeylBook} for the elementary symmetric functions:
\begin{equation}
 \frac{\det D_\Lambda^{(kj)}}{\det D_\Lambda} =  \frac{(-1)^{j-1} E_{k-1} (\Lambda^{\widehat{j}})}{\prod_{l(\not=j)} (1- \Lambda_l / \Lambda_j)},
\end{equation}
where $D_\lambda^{(kj)}$ is $D_\Lambda$ with the $k^\mathrm{th}$ row and $j^\mathrm{th}$ column removed. In this way one immediately retrieves \eqref{mrz-final2} from \eqref{eq:mario}.

The expression \eqref{eq:mario} may be useful for certain applications.
Also, it is easily seen to be directly equivalent to the expression given in \cite{Alf04}.

\subsection{Real case}\label{sec4.2}

\subsubsection{Approach using the eigenvalues of $Q$}

The main trick in deriving the explicit result for $\beta = 2$ was to use the eigen-representation \eqref{mrz-eigen} for the supermatrix $Q$. We are now going to carry out an analogous derivation for $\beta = 1$. The outcome will be somewhat different in that the integrand is no longer meromorphic but has square-root singularities in the radial variables of the boson-boson block.

We use the parametrization
\begin{align}\label{eq:coords}
 Q &= k \begin{bmatrix} Q_0 &0\cr 0 &s \eins_2 \end{bmatrix} k^{-1}, \quad k = \begin{bmatrix} k_0 &0 \cr 0 &\eins_2 \end{bmatrix} \begin{bmatrix} \sqrt{1-\alpha\tilde\alpha} &-\alpha\cr \tilde\alpha ~ &\sqrt{1- \tilde \alpha\alpha} \end{bmatrix} , \\ k_0 &= \begin{bmatrix} \cos\phi &~ -\sin\phi \cr \sin\phi &\cos\phi \end{bmatrix}, \quad \alpha = \begin{bmatrix} \eta & \eta^\ast \cr \xi & \xi^\ast \end{bmatrix} , \quad \tilde\alpha = \begin{bmatrix} \eta^\ast & \xi^\ast \\ -\eta & -\xi \end{bmatrix} , \quad Q_0 = \begin{bmatrix} r_a &0 \cr 0 &r_b \end{bmatrix} \nonumber ,
\end{align}
which has coordinate singularities at $r_a = s$ and $r_b = s$. As before, we will suppress their effects by using the scale invariance of $DQ$ to place the singular locus of the coordinate system on the zero locus of the integrand.

In the coordinates given by \eqref{eq:coords} our invariant integral \eqref{mrz-sb4} takes the form
\begin{align}\label{mrz-4.i}
 \int DQ \, F_1(Q) = \lim_{q\to\infty} &\int \!\!\!\!\!\!\! \int\limits_{ [0,q)^2} \!\! dr_a dr_b \oint\limits_{|s|=q} \!\! \frac{ds} {16\pi \mathrm{i}}\, \frac{|r_a - r_b|} {(r_a - s)^2 (r_b - s)^2} \\ &\times \frac{s}{\sqrt{r_a r_b}} \int\limits_0^{2\pi} \frac{d\phi}{2\pi} \,\partial_\xi \partial_{\xi^\ast} \partial_\eta \partial_{\eta^\ast} (1 + \xi^\ast \xi + \eta^\ast \eta) F_1 \nonumber ,
\end{align}
where $F_1$ is the function
\begin{equation}\label{mrz-f(Q)}
 F_1(Q) = \detg^{n/2}(Q) \exp \left( \frac{\mathrm{i}}{2} \trg XQ \right) \prod_{j=1}^p {\detg}^{-1/2}\left(\eins_4 + \mathrm{i} Q \Lambda_j \right) .
\end{equation}
The expression \eqref{mrz-4.i} holds if $F_1(Q)$ and its derivatives vanish on the singular locus $r_a = s = q \to \infty$ and $r_b = s = q \to \infty$. Observing that our integrand $F_1(Q)$ for $r_{a/b} = s = q$ contains the exponential $\mathrm{e}^{ \mathrm{i} q (\frac{1}{2} x_0 - x_1)}$, we see that $F_1$ has the required property if the parameter range is restricted by $\frac{1}{2} \im\, x_0 > \im\, x_1$. We thus impose this restriction.

It is clear that the integrand \eqref{mrz-f(Q)} does not depend on the angle $\phi$, and all of the dependence on the Grassmann variables is in the factor $\mathrm{e}^{\mathrm{i} \trg X Q /2}$. Hence the result of doing the Grassmann integral is given by
\begin{align}
 &\partial_\xi \partial_{\xi^\ast} \partial_\eta \partial_{\eta^\ast} (1 + \xi^\ast \xi + \eta^\ast \eta) \, \mathrm{e}^{\mathrm{i} \trg XQ / 2} = \mathrm{e}^{\mathrm{i} x_0 (r_a + r_b) / 2 - \mathrm{i} x_1 s} \cr &\times \left( \mathrm{i} (x_0 - x_1) (r_a + r_b - 2 s) - (x_0 - x_1)^2 (r_a - s) (r_b - s) \right) .
\end{align}
A slight reorganization of the integrand then leads to the following expression for the generating function:
\begin{align}\label{mrz-4.ii}
 &Z_1(x_0,x_1) = \left(\frac{x_0}{x_1}\right)^{n-p}\!\!\!\! \lim_{q\to \infty} \int \!\!\!\!\!\!\! \int\limits_{ [0,q)^2} \!\!\! dr_a dr_b
 \,|r_a - r_b| \oint\limits_{|s|=q}\!\! \frac{ds}{16\pi\mathrm{i}}\, \left( \frac{ \sqrt{r_a r_b}}{s} \right)^{n-1}\\ &\times \bigg( \frac{\mathrm{i}(x_0 - x_1) (r_a + r_b -2s)}{(r_a -s)^2(r_b - s)^2} - \frac{(x_0 - x_1)^2}{(r_a-s)(r_b-s)} \bigg) \frac{g_\Lambda(x_1;s)} {\sqrt{g_\Lambda(x_0;r_a) g_\Lambda(x_0;r_b)}} , \nonumber
\end{align}
where $g_\Lambda$ was defined in \eqref{def-gL}.

As before, it is possible to carry out the complex contour integral over the variable $s$. By making partial fraction decompositions to express the integrand in \eqref{mrz-4.ii} by isolated poles $(r_a - s)^{-\ell}$ and $(r_b - s)^{-\ell}$ of degrees $\ell = 1, 2$, and then applying the identity \eqref{eq:cauchy}, we obtain
\begin{align}\label{mrz-4.iii}
 Z_1(x_0,x_1) = \frac{1}{8} \left(\frac{x_0}{x_1}\right)^{n-p} &\int \!\!\!\!\! \int\limits_{\mathbb{R}_+^2} dr_a dr_b \, \frac{|r_a - r_b| \sqrt{r_a r_b}^{\,n-1}} {\sqrt{g_\Lambda(x_0;r_a) g_\Lambda(x_0;r_b)}} \cr &\times \sum_{\nu = 1}^2 (x_0 - x_1)^\nu F^{(\nu)}_\Lambda(x_1; r_a,r_b),
\end{align}
where we have taken the upper limits of the $r$-integrals to infinity, and
\begin{align}
 &F^{(\nu)}_\Lambda(x;r_a,r_b) = \frac{G^{(\nu)}_\Lambda(x;r_a) - G^{(\nu)}_\Lambda(x;r_b)}{r_a - r_b} , \label{Xpress-Fnu} \\
 &G^{(1)}_\Lambda(x;r) = \mathrm{i} \frac{\partial}{\partial r} G^{(2)}_\Lambda (x;r) , \quad G^{(2)}_\Lambda(x;r) = - \frac{g_\Lambda^{[\geq n-1]}(x;r)} {r^{n-1}} . \label{Xpress-Gnu}
\end{align}
The result \eqref{mrz-4.iii} is the analog for $\beta = 1$ of the earlier formula \eqref{mrz-3.vx} for $\beta = 2$.

By recalling the formula \eqref{2.5} we can again deduce an expression for the one-point function, $S_1(x)$. (This requires a process of analytic continuation to remove the restriction $\im\, x_0 > 2\, \im\, x_1$.) The result involves only the $\nu = 1$ term from \eqref{mrz-4.iii} and holds for all $n > p+3$. Unfortunately, we do not know how to do the final two integrals over $r_a$ and $r_b$ in closed form, and an attempt to compute these integrals numerically met with the following complications.

The first complication is that we do not know how to pass from the expression (\ref{mrz-4.iii}--\ref{Xpress-Gnu}) involving $g^>$ to an analog of Eq.\ \eqref{mrz-40} involving $g^<$. Working directly with $g^>$ is numerically expensive because many terms in the Taylor series have to be summed when $r_a$ or $r_b$ are large. On the other hand, if we use the identity $g^> = g - g^<$ then we incur cancelations, with ensuing rounding errors, due to taking the difference of two large numbers.

The second difficulty is that the exponential part of the integrand oscillates. For the most part, these oscillations can be cured by rotating the integration contours for $r_a$, $r_b$ to the positive imaginary axis. Ultimately, however, these variables still have to go to infinity in the direction of the real axis to retain convergence of the integral. The contribution from this ultimate part of the integration contour converges rather slowly.

These problems notwithstanding, we are still able to verify our formula for $S_1(x)$. However, by numerical integration using Wolfram-Mathematical we are not able to produce stable results in a wide range of parameters.
We therefore refrain from even writing down that formula, and abandon now the eigen-representation \eqref{eq:coords}. Instead, we pursue another approach, exploiting the original coordinate system for $Q$ described in Secs.\ \ref{sec3.3} and \ref{sec3.4}.

\subsubsection{Direct approach}

Below, we give two further expressions for the generating function. Their logical order of presentation depends on which of our two formalisms is used.

Adopting the generalized Hubbard-Stratonovich transformation, we start from the integral representation \eqref{3.13} and evaluate the delta-distributions featured in (\ref{3.10a-mrz}--\ref{A.11}). We then carry out the Grassmann integrals (or rather, derivatives) according to the flat measure specified in \eqref{eq3.24} to obtain
\begin{align}\label{eq:56}
 Z_1(x_0,x_1) = \frac{1}{8} \int\!\!\!\!\! \int\limits_{\mathbb{R}_+^2} \! &dr_a dr_b\,\frac{|r_a-r_b|\, \sqrt{r_a r_b}^{\,n-3} \mathrm{e}^{ -(r_a + r_b)/2}}{\prod_{i=1}^p \sqrt{(x_0 -\Lambda_i r_a)(x_0 - \Lambda_i r_b)}} \cr \times \sum_{k=0}^p \frac{(-1)^k x_1^{p-k}}{(n-k)!} &\Bigg( n (n-1) E_k(\Lambda) + \sum_{j\not= l}  \frac{ \Lambda_j^2 \Lambda_l^2 E_{k-2}(\Lambda^{ \widehat{jl}})\,r_a r_b}{(x_0 - \Lambda_j r_a)(x_0 - \Lambda_l r_b)} \\ &+ (n-1) \sum_{j=1}^p \Lambda_j^2 E_{k-1}(\Lambda^{ \widehat{j}}) \left(\frac{r_a} {x_0 - \Lambda_j r_a} + \frac{r_b}{x_0 - \Lambda_j r_b} \right) \Bigg) , \nonumber
\end{align}
where $E_k$, defined in \eqref{eq:ESF}, is understood to vanish when $k < 0$. The notation $\Lambda^{\widehat{jl}}$ means that both $\Lambda_j$ and $\Lambda_l$ are removed from the set $\{\Lambda_i\}_{i = 1,\ldots, p}$.

Alternatively, we can take the formula \eqref{mrz-sb4} from superbosonization and make the substitution $Q \mapsto \mathrm{i} Q / X$. Then, by carrying out the Grassmann integrals in the original parametrization specified in Sec.\ \ref{sec3.4} we obtain
\begin{align}\label{eq:57}
 &Z_1(x_0,x_1) = \frac{(x_1/x_0)^p}{16\pi\mathrm{i}} \int \!\!\!\!\! \int\limits_{\mathbb{R}_+^2}\!\! dr_a dr_b\,\frac{|r_a - r_b|\, \sqrt{r_a r_b}^{\,n-3}}{\sqrt{g_\Lambda (x_0;\mathrm{i}r_a / x_0) g_\Lambda(x_0; \mathrm{i} r_b / x_0)}} \cr &\times\oint\limits_{|s| = 1} \!\!\! ds \; \frac{g_\Lambda(x_1;\mathrm{i}s/x_1)}{s^{n+1}} \Bigg( \sum_{j \not= l} \frac{(\Lambda_j r_a) (\Lambda_j s) (\Lambda_l r_b) (\Lambda_l s)}{(x_0 - \Lambda_j r_a)(x_1 - \Lambda_j s)(x_0 - \Lambda_l r_b)(x_1 - \Lambda_l s)} \cr &\hspace{0.5cm} - (n-1)\sum_{j=1}^p \frac{\Lambda_j s}{x_1 - \Lambda_j s} \left( \frac{ \Lambda_j r_a}{x_0 - \Lambda_j r_a} + \frac{\Lambda_j r_b}{x_0 - \Lambda_j r_b} \right) + n (n-1) \Bigg).
\end{align}
If we go on to carry out the $s$-integral, we again arrive at \eqref{eq:56} above. Note that \eqref{eq:56} and \eqref{eq:57} make immediate sense as convergent integrals for all $x_1 \in \mathbb{C}$ and $x_0 \in \mathbb{C} \setminus \mathbb{R}_+$. A non-trivial check is $Z_1(x_0,x_0) = 1$.

From the expression \eqref{eq:57} it is easy to verify the large-$n$ behavior predicted by the central limit theorem (CLT). Indeed, for large $n$ our integrand develops a saddle point at $r_a = r_b = s \simeq n$ due to the presence of the factor $(\sqrt{r_a r_b} / s)^n \mathrm{e}^{-(r_a + r_b)/2 + s}$. By steepest descent evaluation of the integral around this saddle point one finds
\begin{equation}
 \lim_{n\to\infty} Z_1(n x_0 , n x_1) = \prod_{l=1}^p \frac{x_1 - \Lambda_l}{x_0 - \Lambda_l} ,
\end{equation}
which is the result expected from CLT for the random matrix $(W W^T)_{ij} = \sum_{k=1}^n W_{ik} W_{jk}$ with independent identically distributed matrix elements $W_{jk}$.

\section{Numerical computation of $S_1(x)$} \label{sec.5}

On recalling the basic identity \eqref{2.5}, one immediately produces an expression for the one-point function $S_1(x)$; see below. The result looks more involved than the one obtained from the eigen-representation of $Q$, and the reciprocal square roots visible in Eqs.\ (\ref{eq:56}, \ref{eq:57}) remain serious obstacles for further analytical calculations. From a numerical perspective, however, the new formula has two clear advantages over the previous one: its integrand is free of oscillatory exponentials, and high powers of the integration variables $r_a$ and $r_b$ occur only through the factor $\mathrm{e}^{-(r_a + r_b)/2}$, which is numerically inexpensive to compute. Thus numerical evaluation of the integral \eqref{eq:S1(x)} for $S_1(x)$ is straightforward, provided that the singularities at $r_a, r_b = x / \Lambda_j$ ($j=1,\ldots,p$) are treated with care. We now discuss how to organize this computation.

To prepare our final formula for $S_1(x)$, we observe that the discontinuity $Z(x-\mathrm{i} \varepsilon,x_1) - Z(x+\mathrm{i} \varepsilon, x_1)$ for $\varepsilon \to 0$ arises solely from the product of factors
\begin{equation}\label{eq:sqrts}
 1 / \sqrt{(x - \Lambda_i r_a)(x - \Lambda_i r_b)} ,
\end{equation}
since all other terms in \eqref{eq:56} are single-valued across the real $x$-axis. By the derivation of the result \eqref{eq:56}, the reciprocal square root \eqref{eq:sqrts} has the real value $(x - \Lambda_i r)^{-1}$ for $r_a = r_b = r$. Its values away from the diagonal $r_a = r_b$ are defined by the process of analytic continuation.

For definiteness, let the $\Lambda$-values be ordered by
\begin{displaymath}
 0 < \Lambda_1 < \Lambda_2 < \ldots < \Lambda_p \,.
\end{displaymath}
Fixing any $x > 0$, we partition (a dense open subset of) the domain $\mathbb{R}_+$ for $r_a$ and $r_b$ into $p+1$ intervals $U_0(x) := (0,x/\Lambda_p)$, $U_p(x) := (x/\Lambda_1,\infty)$, and
\begin{equation}
 U_{p-l}(x) := (x/\Lambda_{l+1}, x/\Lambda_l) \qquad (l=1, \ldots, p-1).
\end{equation}
The product over $i = 1, \ldots, p$ of the square roots in \eqref{eq:sqrts} is then discontinuous across the $x$-axis only for $(r_a,r_b) \in U_l(x) \times U_{l^\prime}(x)$ with $l+l^\prime$ an odd number. Thus we arrive at an integral representation for $S_1(x)$ of the form
\begin{align}\label{eq:S1(x)}
 S_1(x) &= \!\! \sum_{l+l^\prime \text{odd}} \, \int\limits_{U_l(x)} \!\!\!\! dr_a \int\limits_{U_{l^\prime}(x)}\!\!\!\!\! dr_b\; \frac{(-1)^{(l+l^\prime - 1)/2}} {\prod_{i=1}^p \sqrt{|(x -\Lambda_i r_a)(x - \Lambda_i r_b)|}} \\ &\times \bigg( f_{x,\Lambda} (r_a,r_b) + \sum_{j=1}^p \frac{f_{x ,\Lambda}^{(j)}(r_a,r_b)} {x - \Lambda_j r_a} + \sum_{j\not= j^\prime}  \frac{f_{x,\Lambda}^{(j j^\prime)}(r_a,r_b)}{(x - \Lambda_j r_a)(x - \Lambda_{j^\prime} r_b)} \bigg) \nonumber .
\end{align}
The outer sum is over the pairs $(l,l^\prime) \in \{0,1,\ldots,p\}^2$ with odd sum $l + l^\prime$, and the $f_{x,\Lambda}^{(..)}(r_a,r_b)$ are certain real-valued analytic functions, which do not need to be written down as they are easy to read off from \eqref{eq:56} using \eqref{2.5}.

What we do need is some further discussion of the precise meaning of the integral \eqref{eq:S1(x)}. There is no problem with the first summand (containing $f_{x,\Lambda}$), as the square root singularities from \eqref{eq:sqrts} are integrable. However, in the second and third summand of \eqref{eq:S1(x)} the order of the singularity is enhanced from $(x - \Lambda_j r)^{-1/2}$ to $(x - \Lambda_j r)^{-3/2}$ for one ($j$) or two ($j \not= j^\prime$) factors of the product over $i = 1, \ldots, p$. These terms have to be properly understood as principal-value integrals by their definition as a limit $\varepsilon \to 0$. To render them suitable for numerical evaluation, we use a partial integration identity of the sort
\begin{equation}\label{8.8}
 \frac{1}{2} \lim_{\varepsilon\to 0}\int\limits_{-\delta_-}^{\delta_+} \frac{g(r)\, dr} {(r- \mathrm{i}\varepsilon)^{3/2}} = - \frac{g(r)} {r^{1/2}} \Bigg\vert_{r = -\delta_-}^{r = \delta_+} + \int\limits_{ -\delta_-}^{\delta_+} \frac{\partial_r g(r) \, dr} {r^{1/2}} ,
\end{equation}
which holds for any smooth function $g(r)$. We apply this identity to the present situation with $r \equiv r_a - \Lambda_j / x$ and/or $r \equiv r_b - \Lambda_{j^\prime} / x$ (and suitable values for $\delta_\pm$), and thereby regularize the integral around the $3/2$ singularities. In order for this regularization to be well-defined, we need the diagonal $r_a = r_b$ to lie (as it does) outside the domain of integration, so that our integrand $g(r)$ containing the factor $|r_a - r_b|$ is in fact smooth.

\begin{figure*}[ht!]
\centering
\includegraphics[scale=0.25]{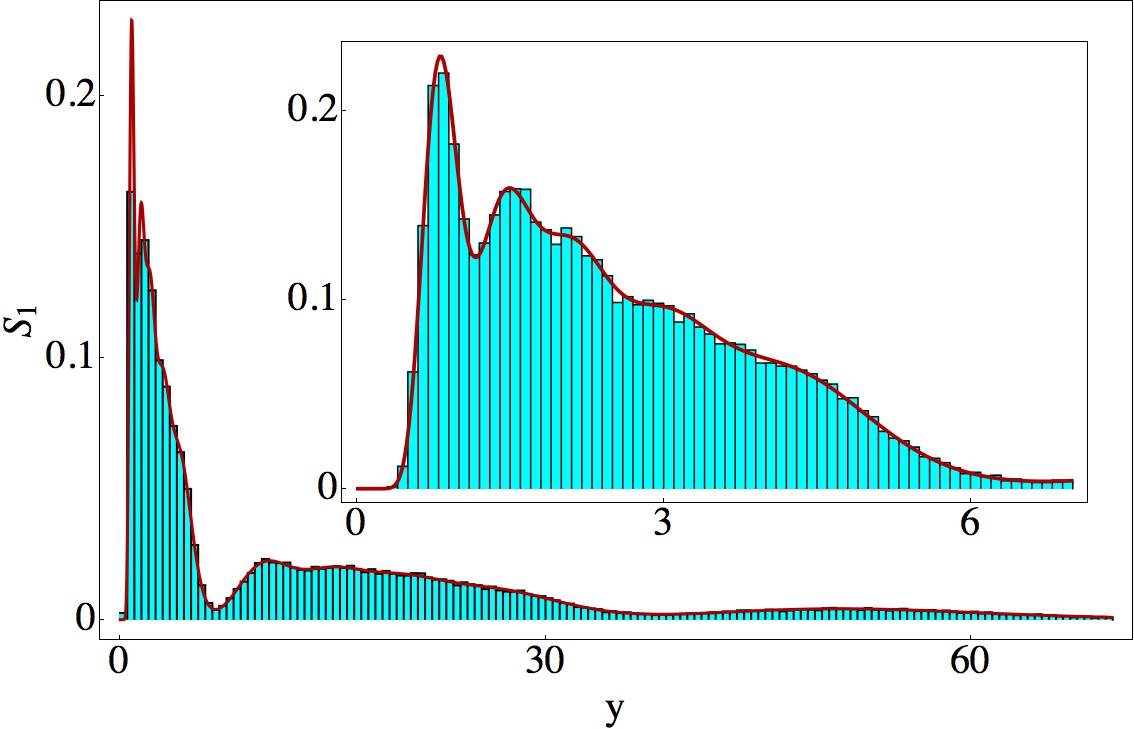}
\caption{One-point function $S_1(x)$ for the parameter values $p=10$ and $n=50$. The set of $\{\Lambda_j\}_{j=1, \ldots, 10}$ is $\{1, 0.49, 0.4225, 0.36, 0.25, 0.09, 0.0729, 0.0529, 0.04, 0.0225\}$. The solid line is the result obtained by numerical integration of the analytical formula \eqref{eq:S1(x)}. The bin size of the histogram is 0.7 for the large figure and 0.1 for the inset. The inset magnifies $S_1(x)$ near the origin.} \label{fig1}
\end{figure*}
\begin{figure*}[ht!]
\centering
\includegraphics[scale=0.25]{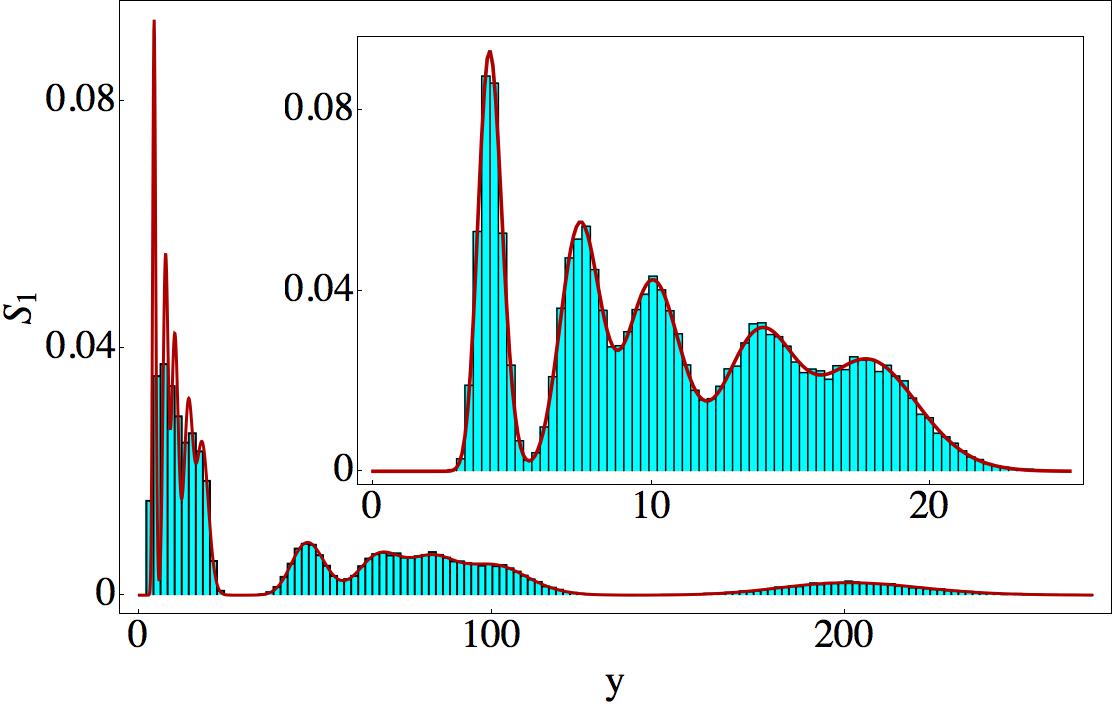}
\caption{$S_1(x)$ for the same values of $p$ and $\{\Lambda_j\}$ but a larger parameter $n = 200$. The solid line is again the result obtained from Eq.\ \eqref{eq:S1(x)}. The bin size is 2 for the large figure and 0.3 for the inset. Note that because of the bigger value of $n$ the peaks are more pronounced (than those in Fig.\ \ref{fig1}) around the CLT-values $n\Lambda_j$.} \label{fig2}
\end{figure*}

By using the commercial software
\textsc{Mathematica}\textsuperscript{\textregistered} \cite{math} we
now numerically evaluate the integrals in Eq.\ \eqref{eq:S1(x)}. In
Figs.\ \ref{fig1} and \ref{fig2} we show examples for $(p,n) = (10,50)$ and $(10,200)$, respectively. Our result (solid line) is compared with a Monte-Carlo simulation (histogram) using an ensemble of $10^5$ random matrices. We see that the agreement between the Monte-Carlo simulation and our result from \eqref{eq:S1(x)} is perfect.

We observe that the last peak, centered around $x = 50$ in Fig.\ \ref{fig1} and $x = 200$ in Fig.\ \ref{fig2}, lies near $n \Lambda_p = n$, in agreement with what is expected by the central limit theorem (CLT).
%
%
The other peaks are slightly shifted from their asymptotic CLT-positions  $n \Lambda_j$. With increasing $n$ these shifts become smaller and the peaks become more pronounced, as it should be.

\section{Summary}\label{sec.6}

We have developed a supersymmetry approach to derive exact expressions for the one-point function of real and complex Wishart correlation matrices. Supersymmetry got us around a difficult group integral which arises in the traditional approach for the real case. The crucial advantage of the supersymmetry approach is the drastic reduction of the number of integrals. For both cases, complex and real, we showed how to express the one-point function as an eigenvalue integral. In the complex case the supermatrix has dimension $2 \times 2$ and thus two eigenvalues. We carried out both integrals and demonstrated the equivalence of our result to that of Ref.\ \cite{Alf04}.

In the more demanding real case ($\beta = 1$), the supermatrix has dimension $4 \times 4$ and 3 distinct eigenvalues. While the Efetov-Wegner boundary terms due to diagonalization for $\beta = 1$ have never been given in explicit form, we showed how to suppress them by a variable substitution that pushes the coordinate singularities outside the domain of integration. One of the three eigenvalue integrals is easily done by residue calculus. We do not know how to calculate the remaining twofold integral by analytical means, and an attempt to compute it numerically met with some complications. We therefore abandoned the coordinate system given by the eigenvalues and turned to a direct approach using standard coordinates. Thus we produced a second formula for $S_1(x)$, still as a twofold integral, which we were able to compute numerically in a stable and efficient way. We also illustrated our result by comparison with a Monte-Carlo simulation.

Previous approaches to the real case had resulted in slowly converging series of Jack or zonal polynomials. We believe that our formula represents a considerable improvement over these previous results. From a conceptual viewpoint, one might say that our result re-sums a multiple infinite series of Jack polynomials (and integrals thereof) in a non-trivial real case.

Thus we hope that we have demonstrated that the supersymmetry method is a powerful tool to tackle problems in multivariate statistics.

\section*{Acknowledgements}

We thank R.\ Sprik for fruitful discussions on the theoretical and
experimental aspects of the one-point function. We are also grateful to
A.\ Hucht, H.\ Kohler and R.\ Sch\"afer for helpful comments. One of us (TG) thanks the organizers and participants of the Program in 2006 on ``High Dimensional Inference and Random Matrices'' at SAMSI, Research Triangle Park, North Carolina (USA), for drawing his attention to this problem. We acknowledge support from the German Research Council (DFG) via the Sonderforschungsbereich-Transregio 12 ``Symmetries and Universality in Mesoscopic Systems''.


\begin{thebibliography}{10}

\bibitem{GMW98}
T.~Guhr, A.~M{\"u}ller-Groeling, and H.A. Weidenm{\"u}ller.
\newblock {\em Phys. Rep.}, 299:189, 1998.

\bibitem{Ver04}
A.~Tulino and S.~Verd\'u.
\newblock Found. trends commun. inf. theory 1.
\newblock 1, 2004.

\bibitem{Haa01}
F.~Haake.
\newblock {\em Quantum Signatures of Chaos}.
\newblock Springer Verlag, Berlin, 2nd edition, 2004.

\bibitem{Voi05}
J.~Voit.
\newblock {\em The statistical mechanics of financial markets}.
\newblock Springer Verlag, Berlin, 3rd edition, 2005.

\bibitem{lal99}
L.~Laurent, P.~Cizeau, J.P.~Bouchaud, and M.~Potters.
\newblock {\em Phys. Rev. Lett.}, 83:1467, 1999.

\bibitem{Spr08}
R.~Sprik, A.~Tourin, J.~de~Rosny, and M.~Fink.
\newblock {\em Phys. Rev. B}, 78:012202, 2008.

\bibitem{Mui82}
R.J.~Muirhead.
\newblock {\em Aspects of Multivariate Statistical Theory}.
\newblock Wiley, New York, 1st edition, 1982.

\bibitem{Alf04}
G.~Alfano, A.~Tulino, A.~Lozano, and S.~Verd\'u.
\newblock {\em Proc. IEEE 8. Int. Symp. on Spread Spectrum Tech. and Applications (ISSSTA '04) (IEEE, Bellingham,WA, 2004)}, pp. 515.

\bibitem{Sil95}
J.W.~Silverstein.
\newblock {\em J. Multivariate Anal.}, 55:331, 1995.

\bibitem{VP2010}
Vinayak and Akhilesh Pandey.
\newblock {\em Phys. Rev. E}, 81(3):036202, 2010.

\bibitem{Efe97}
K.B.~Efetov.
\newblock {\em Supersymmetry in Disorder and Chaos}.
\newblock Cambridge University Press, Cambridge, 1st edition, 1997.

\bibitem{Verb04}
J.J.M.~Verbaarschot.
\newblock {\em AIP Conf. Proc}, 744:277, 2004.

\bibitem{RKG10/1}
C.~Recher, M.~Kieburg, and T.~Guhr.
\newblock {\em Phys. Rev. Lett.}, 105(24):244101, 2010.

\bibitem{Guh06}
T.~Guhr.
\newblock {\em J. Phys.}, A 39:13191, 2006.

\bibitem{KGG08}
M.~Kieburg, J.~Gr{\"o}nqvist, and T.~Guhr.
\newblock {\em J. Phys. A}, 42:275205, 2009.

\bibitem{Som07}
H.-J.~Sommers.
\newblock {\em Acta Phys. Pol.}, B 38:1001, 2007.

\bibitem{LSZ07}
P.~Littelmann, H.-J.~Sommers, and M.R.~Zirnbauer.
\newblock {\em Commun. Math. Phys.}, 283:343, 2008.

\bibitem{KSG08}
M.~Kieburg, H.J.~Sommers, and T.~Guhr.
\newblock {\em J. Phys. A}, 42:275206, 2009.

\bibitem{Ver94/2}
J.J.M.~Verbaarschot.
\newblock {\em Phys. Rev. Lett.}, 72:2531, 1994.

\bibitem{Har58}
Harish-Chandra.
\newblock {\em Am. J. Math.}, 80:241, 1958.

\bibitem{ItzZub80}
C.~Itzykson and J.B.~Zuber.
\newblock {\em J. Math. Phys.}, 21:411, 1980.

\bibitem{EST04}
K.B.~Efetov, G.~Schwiete, and K.~Takahashi.
\newblock {\em Phys. Rev. Lett.}, 92:026807, 2004.

\bibitem{VWZ85}
J.J.M.~Verbaarschot, H.A.~Weidenm{\"u}ller, and M.R.~Zirnbauer.
\newblock {\em Phys. Rep.}, 129:367, 1985.

\bibitem{Guh10}
T.~Guhr.
\newblock {\em  The Oxford Handbook of Random Matrix Theory, G. Akemann, J. Baik and P. Di Francesco (eds.), Chapter 7: Supersymmetry,  pages 135 - 154}.
\newblock Oxford University Press, Oxford, 2011, preprint: arXiv:1005.0979.

\bibitem{Ber87}
F.A.~Berezin.
\newblock {\em Introduction to Superanalysis}.
\newblock D. Reidel Publishing Company, Dordrecht, 1987.

\bibitem{Zir96}
M.R.~Zirnbauer.
\newblock {\em J. Math. Phys.}, 37:4986, 1996.



\bibitem{Fyo02}
Y.V.~Fyodorov.
\newblock {\em Nucl. Phys.}, B 621:643, 2002.






\bibitem{math}
Wolfram~Research Inc.
\newblock \textsc{Mathematica}~version 7.0.
\newblock 2008.


\bibitem{WeylBook} H.~Weyl.
\newblock {\em The Classical Groups: Their Invariants and Representations}.
\newblock Princeton University Press, Princeton, 1939.

\end{thebibliography}
\end{document}